\documentclass[leqno]{article}

\usepackage{graphicx}
\usepackage{amsmath}
\usepackage{amssymb}
\usepackage{amsfonts}
\usepackage{aeguill}
\usepackage{float}
\usepackage[english]{babel}

\newcommand{\ds}{\displaystyle}
\newcommand{\ub}{{\bar  u}}
\newcommand{\vb}{{\bar  v}}
\newcommand{\rhob}{{\bar  \rho}}
\newcommand{\cb}{{\bar  c}}

\newcommand{\eulerhathat}{{\hat {\cal P}}}

\renewcommand{\vec}[1]{\mbox{\boldmath $#1$}}  
\newcommand{\R}{\mathbb{R}}

\newtheorem{proposition}{Proposition}  
\newtheorem{theorem}{Theorem}

\newtheorem{algorithm}{ALGORITHM}  
  
\newcommand{\QED}{\hspace*{\fill}\rule{2.5mm}{2.5mm}}  
\newenvironment{proof}{{\bf Proof\ }}{\QED\\}  
  
\begin{document}

\title{A New Domain Decomposition Method for the Compressible Euler Equations}  
  
\author{   
V.~Dolean\footnote{Univ. d'\'Evry and CMAP, Ecole Polytechnique, 91128  
Palaiseau, France,  dolean@cmap.polytechnique.fr},~  
F.~Nataf\footnote{CMAP, CNRS UMR 7641, Ecole Polytechnique, 91128 Palaiseau Cedex, France, nataf@cmap.polytechnique.fr}
}  
  
\renewcommand{\theequation}{\arabic{section}.\arabic{equation}}  
\makeatletter  
\@addtoreset{equation}{section}  
\makeatother

\maketitle  

\begin{abstract}  
In this work we design a new domain decomposition method for the Euler equations in $2$ dimensions. The basis is the equivalence via the Smith factorization with a third order scalar equation to whom we can apply an algorithm inspired from the Robin-Robin preconditioner for the convection-diffusion equation. Afterwards we translate it into an algorithm for the initial system and prove that at the continuous level and for a decomposition into $2$ sub-domains, it converges in $2$ iterations. This property cannot be preserved strictly at discrete level and for arbitrary domain decompositions but we still have numerical results which confirm a very good stability with respect to the various parameters of the problem (mesh size, Mach number, $\ldots$). 
\end{abstract}  

\tableofcontents  
\bibliographystyle{alpha}

\section{Introduction}\label{intro}
The need of using domain decomposition methods when solving  partial differential equations is nowadays more and more obvious. The challenge is now the acceleration of these methods. Different possibilities were studied such as the use of optimized interface conditions on the artificial boundaries between subdomains or the preconditioning of a substructured system defined at the interface. The former were widely studied and analyzed for scalar problems such as elliptic equations in \cite{Lions:1990:SAM,Enquist:1998:ABC}, for the Helmholtz equation in \cite{Benamou:1997:DDH,Chevalier:1998:SMO,Gander:2001:OSH,Lee:2005:NOD}  convection-diffusion problems in \cite{Japhet:2000:OO2}. For time dependent problems and local times steps, see for instance \cite{Gander:2001:OSW,Gander:2003:MRO}. The preconditioning methods have also known a wide developpement in the last decade. The Neumann-Neumann algorithms for symmetric second order problems \cite{Bourgat:1989:VFA,Mandel:92:BDD,DeRoeck:91:DDP} has been the subject of numerous works, see \cite{toselli-widlund:04} and references therein. An extension of these algorithms to non-symmetric scalar problems (the so called Robin-Robin algorithms) has been done in \cite{Achdou:2000:DDP,Giorda-etal:DDP:04} for advection-diffusion problems.\\

The generalization of these kinds of methods to the system of the Euler equations is difficult in the subsonic case in dimensions equal or higher to two. As far as Schwarz algorithms are concerned, when dealing   with supersonic  flows, whatever  the   space  dimension is, imposing   the appropriate characteristic  variables as interface conditions leads to a convergence of the  algorithm which is  optimal with regards to the number of subdomains.   This  property is generally lost  for subsonic flows except for the case of one-dimensional  problems, where  the optimality is  expressed by the fact that the number of iterations is equal to the  number of  subdomains  (see  Bj{\o}rhus \cite{Bjorhus:1995:NCD} and Quarteroni \cite{qua:90} for   more details). In the subsonic case and in two or three dimensions, we can find a formulation with classical (natural) transmission conditions in \cite{qua:90,Cai:1998:MORAS,quarteroni-stolcis:96} or with more general interface conditions in \cite{Clerc:98:NSM} and optimized transmission conditions in \cite{Dolean:2004:OSM}. The analysis of such algorithms applied to systems proved to be very different from  the scalar case, see \cite{Dolean-etal:02:SYS2,Dolean-etal:04:SYS}. As far as preconditioning methods are concerned, to our knowledge, no extension to the Euler equations was done. \\ 

In this paper, we consider a preconditioning technique for the system of the compressible Euler equations in the subsonic case. The paper is organized as follows: in Section \ref{equiv} we will first show the equivalence between the 2D Euler equations and a third order scalar problem, which is quite natural by considering a Smith factorization of this system, see \cite{Wloka:1995:BVP} or \cite{Gantmacher:1966:TM}.  In Section \ref{schwarz-scalar} we define an optimal algorithm for the third order scalar equation. It is inspired from the idea of the Robin-Robin algorithm \cite{Achdou:2000:DDP} applied to a convection-diffusion problem. We also prove by using a Fourier analysis that this algorithm converges in two iterations. Afterwards in Section \ref{schwarz-system} we back-transform it and define  the corresponding algorithm applied to the Euler system. All the previous results have been obtained at the continuous level and for a decomposition into $2$ unbounded subdomains. After a discretization in a bounded domain we cannot expect that these properties to be preserved exactly. Still we can show in Section \ref{discrete} by a discrete convergence analysis that the expected results should be very good. In Section \ref{results}, numerical results confirm the very good stability of the algorithm with respect to the various parameters of the problem (mesh size, Mach number, $\ldots$). 
 
\section{A third order scalar problem}\label{equiv}
In this section we will show the equivalence between the linearized Euler system and a third order scalar equation. The motivation for this transformation is that a new algorithm is easier to design for a scalar equation than for a system of partial differential equations. 
\subsection{The compressible 2D Euler equations}
In the following we will concentrate ourselves on the conservative Euler equations in two-dimensions:

\begin{equation}
\label{euler1}
\displaystyle\frac{\partial W}{\partial t} +
{\nabla}.\vec F(W) = 0 
~~,~~ 
W = \left (\rho,\rho\vec V,E \right )^t.
\end{equation}
In the above expressions, $\rho$ is the density, $\vec V=(u,v)^t$ is
the velocity vector,  $E$ is the total  energy per unit  of volume and
$p$ is the pressure. In  equation  (\ref{euler1}), $W = W(\bf   x ,  t)$  is the  vector of conservative variables, $\bf x$ and  $t$ respectively denote the space and  time variables and   $\vec F(W)=\left (F_1(W),F_2(W) \right
)^T$ is the conservative flux vector whose components are given by

$$
F_1(W) = 
\left (\rho u, \rho u{^2} + p,\rho uv,u(E+p) \right )^t ,~~~
F_2(W) = 
\left (\rho v, \rho uv, \rho v{^2} + p, v(E+p)\right )^t.
$$
The pressure is deduced from the other variables
using the state equation for a perfect gas $p = (\gamma_s - 1)(E - \frac{1}{2}\rho\parallel\vec V\parallel^{2})$ where $\gamma_s$ is the ratio of the specific heats ($\gamma_s = 1.4$ for the air).

\subsection{Equivalence of the Euler system to a scalar equation}
The starting point of our  analysis is given by the linearized form of the
Euler equations (\ref{euler1}) written in primitive variables $(p,u,v,S)$. In the following we suppose that the flow is isentropic, which allows us to drop the equation of the entropy (which is totally decoupled with respect to the others). We denote by $W=(P,U,V)^T$ the vector of unknowns and by $A$ and $B$ the jacobian matrices of the fluxes $F_i(w)$ to whom we already applied the variable change from conservative to primitive variables. In the following, we shall denote by $\bar c$ the speed of the sound and we consider the linearized form (we will mark by the bar symbol, the constant state around which we linearize) of the Euler equations:
 
\begin{equation}\label{euler2}
{\cal P}W\equiv \frac{W}{\Delta t} + A\partial_xW+B\partial_yW=f
\end{equation}
characterized by the following jacobian matrices:
\begin{equation}\label{a12}
\begin{array}{cc}
A =\left ( 
\begin{array}{ccc}
\bar u & \bar\rho \bar c^2 & 0  \\ [2ex]
1/\bar\rho & \bar u & 0 \\ [2ex]
0 & 0 & \bar u 
\end{array}
\right ) 
& 
\hspace{1.0cm} 
B = 
\left ( 
\begin{array}{ccc}
\bar v & 0 & \bar\rho \bar c^2   \\ [2ex]
0 & \bar v & 0 \\ [2ex]
1/\bar\rho & 0 & \bar v 
\end{array}
\right ) \\
\end{array}
\end{equation}
We can re-write the system (\ref{euler2}) by denoting $\beta=\frac{1}{\Delta t}>0$ under the form
\begin{equation}\label{euler3}
{\cal P} W\equiv \left(\beta I+A\partial_x+B\partial_y\right)W =f
\end{equation}
We will study this system with the help of the Smith factorization. 
\subsubsection{Smith factorization}
We first recall the definition of the Smith factorization of a matrix with polynomial entries and apply it to systems of PDEs, see \cite{Gantmacher:1966:TM2,Gantmacher:1966:TM1,Gantmacher:1998:TM} or \cite{Wloka:1995:BVP} and references therein. \\
\begin{theorem}
Let $n$ be an integer and $C$ an invertible $n\times n$ matrix with polynomial entries in the variable $\lambda$ : $C=(c_{ij}(\lambda))_{1\le i,j\le n}$. \\
Then, there exist three matrices with polynomial entries $E$, $D$ and $F$ with the following properties:
\begin{itemize}
\item det($E$)=det($F$)=1
\item $D$ is a diagonal matrix.
\item $C=EDF$.
\end{itemize}
Moreover, $D$ is uniquely defined up to a reordering and multiplication of each entry by a constant by  a formula defined as follows. Let $1\le k\le n$,
\begin{itemize}
\item  $S_k$ is the set of all the submatrices of order $k\times k$ extracted from $C$.
\item $Det_k=\{\text{Det}(B_k)\backslash B_k\in S_k  \}$
\item $LD_k$ is the largest common divisor of the set of polynomials $Det_k$. 
\end{itemize}
Then,
\begin{equation}
D_{kk}(\lambda)=\frac{LD_k(\lambda)}{LD_{k-1}(\lambda)}, \ \  1\le k\le n
\end{equation}
(by convention, $LD_0=1$). 
\end{theorem}

\paragraph{Application to the Euler system}
We first take formally the Fourier transform of the system (\ref{euler3}) with respect to $y$ (the dual variable is $\xi$). We keep the partial derivatives in $x$ since in the sequel we shall consider a domain decomposition with an interface whose normal is in the $x$ direction. We note
\begin{equation}\label{eq:eulerhathat}
\eulerhathat = 
\left(\begin{array}{ccc} \beta + \ub \partial_x + i\xi \vb  &  \rhob \cb^2 \partial_x &   i \rhob \cb^2 \xi \\
                         \frac{1}{\rhob}\partial_x & \beta + \ub \partial_x + i\xi\vb  & 0 \\
                      \frac{i\xi}{\rhob} & 0 & \beta +\ub \partial_x+ i\vb \xi 
\end{array} \right)
\end{equation}

We can perform a Smith factorization of  $\eulerhathat$ by considering it as a matrix with polynomials in $\partial_x$ entries.  We have  
\begin{equation}\label{eq:smithfactorization}
\eulerhathat=EDF
\end{equation}
 where
\begin{equation}\label{eq:Dsmith}
D=\left(\begin{array}{ccc} 1 & 0 & 0 \\
                           0 & 1 & 0 \\
                           0 & 0 & \ds\hat{{\cal L}}   \hat{{\cal G}}
        \end{array}\right) 
\end{equation}
and 
\[
E= \frac{1}{(\ub(\cb^2-\ub^2))^{1/3}}\left(\begin{array}{ccc} i \rhob\cb^2\xi & 0 & 0 \\
                           0 & \ub & 0 \\
       \beta +\ub\partial_x+i \vb \xi & E_2 &\ds \frac{\cb^2-\ub^2}{i\xi\rhob\cb^2}
        \end{array}\right)
\]
and  
\begin{equation}\label{sec:F}
F= -\left(\begin{array}{ccc} 
\ds \frac{\beta+\ub\partial_x+i\xi\vb}{i\xi\rhob \cb^2} & \ds\frac{\partial_x}{i\xi} & 1\\
\ds \frac{\partial_x}{\rhob\ub} &\ds \frac{\beta+\ub\partial_x+i\xi\vb}{\ub} & 0\\
\ds \frac{\ub}{\beta+i\xi\vb} & \ds\frac{\rhob\ub^2}{\beta+i\xi\vb} & 0
   \end{array}\right)
\end{equation}
where
\[
E_2=\ub\frac{(-\ub\cb^2+\ub^3)\partial_{xx}+(2\ub^2-\cb^2)(\beta+i\xi\vb)\partial_x+\ub((\beta+i\xi\vb)^2+\xi^2\cb^2)}{\cb^2(i\beta+i\xi\vb)},
\]
\begin{equation}
  \label{eq:Grond}
{\hat {\cal G}}=\beta+\ub\partial_x+i\xi\vb
\end{equation}
 and 
\begin{equation}
  \label{eq:Lrond}
  {\hat {\cal L}}= \beta^2+2i\xi \ub\vb\partial_{x}+2\beta(\ub\partial_x+i\xi\vb)+(\cb^2-\vb^2)\xi^2  - (\cb^2-\ub^2)\partial_{xx}
\end{equation}
The operators showing up in the diagonal matrix have a physical meaning:
\[
{\cal G}=\beta +\ub\partial_x+\vb\partial_y
\] 
is a first order transport operator where the time derivative is replaced by $\beta$ and 
\[
{\cal L}= \beta^2+2 \ub\vb\partial_{xy}+2\beta(\ub\partial_x+\vb\partial_y)-(\cb^2-\vb^2)\partial_{yy}  - (\cb^2-\ub^2)\partial_{xx}
\] 
is the advective wave operator where $\partial_t^l$ is replaced by $\beta^l$ for $l=1,2$.\\

 Equation \eqref{eq:Dsmith} suggests that the derivation of a domain decomposition method (DDM) for the third order operator ${\cal L G}$ is a key ingredient for a DDM for the compressible Euler equations. 
\section{A new algorithm applied to a scalar third order problem}\label{schwarz-scalar}
In this section we will describe a new algorithm applied to the third order operator found in section \ref{equiv}. We want to solve
\begin{equation}
  \label{eq:grondl}
  {\cal L} {\cal G}(Q)=g
\end{equation}
where $Q$ is scalar unknown function and $g$ is a given right hand side. The algorithm will be based on the Robin-Robin algorithm \cite{Achdou:1997:RPA,Achdou:2000:DDP} for the convection-diffusion problem. Then we will prove its convergence in $2$ iterations. We first note that the elliptic operator ${\cal L}$ can also be written as:
\begin{equation}\label{eq:Arond}
{\cal L}=-div(A\nabla)+{\bf a}\nabla+\beta^2,\, A=\left(\begin{array}{cc}\bar c^2-\bar u^2 & -\bar u\bar v\\ -\bar u\bar v & \bar c^2-\bar v^2\end{array}\right) \text{ where } {\bf a}=2\beta(\bar u,\bar v)
\end{equation}
Without loss of generality we assume in the sequel that the flow is subsonic and that $\bar u > 0$ and thus we have $0 < \bar u < \bar c$. 
\subsection{The algorithm for a two-domain decomposition}\label{sec:algo-two-subd}
We consider now a decomposition of the plane $\mathbb{R}^2$ into two non-overlapping sub-domains $\Omega_1=(-\infty,0)\times \R$ and $\Omega_2=(0,\infty,0)\times \R$. The interface is $\Gamma=\{x=0\}$. The outward normal to domain $\Omega_i$ is denoted $\bf n_i$, $i=1,2$. Let $Q^{i,k},\,i=1,2$ represent the approximation to the solution in subdomain $i$ at the iteration $k$ of the algorithm. We define the following algorithm:
\begin{algorithm}\label{algo1}
We choose the initial values $Q^{1,0}$ and $Q^{2,0}$ such that ${\cal G}Q^{1,0}={\cal G}Q^{2,0}$ and we compute $(Q^{i,k+1})_{i=1,2}$ from $(Q^{i,k})_{i=1,2}$ by the following iterative procedure:\\
{\bf Correction step} We compute the corrections $\tilde Q^{1,k}$ and $\tilde Q^{2,k}$ as solution of the homogeneous local problems:
\begin{equation}\label{step1}
\begin{array}{cc}
\left\{\begin{array}{l}
{\cal LG}\tilde Q^{1,k}=0\text{ in }\Omega_1,\\[2ex]
(A\nabla-\frac{1}{2}{\bf a}){\cal G}\tilde Q^{1,k}\cdot{\bf n_1}=\gamma^k,\text{ on }\Gamma.
\end{array}\right.
&
\left\{\begin{array}{l}
{\cal LG}\tilde Q^{2,k}=0\text{ in }\Omega_2,\\[2ex]
(A\nabla-\frac{1}{2}{\bf a}){\cal G}\tilde Q^{2,k}\cdot{\bf n_2}=\gamma^k,\text{ on }\Gamma,\\[2ex]
\tilde Q^{2,k}=0,\text{ on }\Gamma.
\end{array}\right.
\end{array}
\end{equation} 
where $\gamma^k=-\frac{1}{2}\left[A\nabla{\cal G}Q^{1,k}\cdot{\bf n_1}+A\nabla {\cal G}Q^{2,k}\cdot{\bf n_2}\right]$.\\
{\bf Update step}.We update $Q^{1,k+1}$ and $Q^{2,k+1}$ by solving the local problems:
\begin{equation}\label{step2}
\begin{array}{cc}
\left\{\begin{array}{l}
{\cal LG}Q^{1,k+1}=g,\text{ in }\Omega_1,\\[2ex]
{\cal G}Q^{1,k+1}={\cal G}Q^{1,k}+\delta^k,\text{ on }\Gamma.
\end{array}\right.
&
\left\{\begin{array}{l}
{\cal LG}\tilde Q^{2,k+1}=g,\text{ in }\Omega_2,\\[2ex]
{\cal G}Q^{2,k+1}={\cal G}Q^{2,k}+\delta^k,\text{ on }\Gamma,\\[2ex]
Q^{2,k+1}=Q^{1,k}+\tilde Q^{1,k},\text{ on }\Gamma.
\end{array}\right.
\end{array}
\end{equation} 
where $\delta^k=\frac{1}{2}\left[{\cal G}\tilde Q^{1,k}+{\cal G}\tilde Q^{2,k}\right]$.
\end{algorithm}

\begin{proposition}
Algorithm \ref{algo1} converges in $2$ iterations.
\end{proposition}
\begin{proof} We use the Fourier transform technique. For the sake of the analysis we consider the previous algorithm written in term of error vector $e^{i,k}(x,y)=(Q^{i,k}-Q|_{\Omega_i})(x,y)$. The error  $(e^{i,k})_{i=1,2}$ satisfies Algorithm~\ref{algo1} with $g=0$.\\
We will first describe what happens locally inside each subdomain after proceeding to a Fourier transform in the $y$ direction and then we prove the convergence of the algorithm \ref{algo1} by computing in Fourier space the effect of the correction and the update steps. We denote by $\hat e(x,\xi)$  the Fourier transform of a function $e(x,y)$:
$$
\hat e(x,\xi) = 
\displaystyle\int_{\mathbb{R}}e(x,y) e^{-i\xi y}dy
$$
We first study solutions to the homogeneous equation ${\cal L}{\cal G}(e_i)=0$ in domain $\Omega_i$, $i=1,2$. We take its Fourier in the $y$ direction and get:
\begin{equation}\label{eq1}
\hat{\cal L}\hat{\cal G}\hat e^i=(\beta+\bar u\partial_x+i\xi \bar v)\,((\bar u^2-\bar c^2)\partial_{xx}+2\bar u(\beta+i\xi \bar v)\partial_x+\xi^2+(\beta+i\xi \bar v)^2)\hat e^i=0
\end{equation}
We seek the solution in the form $\hat e^i(x,\xi)=e^{\lambda(\xi)x}$.and we find three possible values for $\lambda$:
\begin{equation}
\lambda_{1,2}(\xi)=\frac{\bar u(\beta+i\xi \bar v)\pm \bar c\sqrt{(\beta+i\xi \bar v)^2+\xi^2(\bar c^2-\bar u^2)})}{\bar c^2-\bar u^2},\,\lambda_3(\xi)=-\frac{\beta+i\xi \bar v}{\bar u}.
\end{equation}
therefore the solution writes $\hat e^i(x,\xi)=\alpha_{1i}(\xi)e^{\lambda_1(\xi)x}+\alpha_{2i}(\xi)e^{\lambda_2(\xi)x}+\alpha_{3i}(\xi)e^{\lambda_3(\xi)x}$. We also impose that $\hat e_1$ (resp. $\hat e_2$) is bounded as $x$ tends to $-\infty$ (resp. $\infty$). Taking into account the sign of the real parts of $(\lambda_j)_{j=1,2,3}$ it means that we have:
\begin{equation}\label{sols}
\hat e^1(x,\xi)=\alpha_1e^{\lambda_1(\xi)x}\text{ and }\hat e^2(x,\xi)=\alpha_2e^{\lambda_2(\xi)x}+\alpha_3e^{\lambda_3(\xi)x}
\end{equation}
In order to ease the notations we call:
\begin{equation}
a(\xi) = \beta+i\xi \bar v \text{ and }R(\xi)=\sqrt{(\beta+i\xi \bar v)^2+\xi^2(\bar c^2-\bar u^2)}
\end{equation}
The initial guesses of the algorithm does not satisfy a specific partial differential equation. Therefore, it is possible to use formula \eqref{sols} only for $k\ge 1$. With obvious notations, we write:
\begin{equation}\label{eq:e}
\hat e^{1,k}(x,\xi)=\alpha_1^k e^{\lambda_1(\xi)x}\text{ and }\hat e^{2,k}(x,\xi)=\alpha_2^k e^{\lambda_2(\xi)x}+\alpha_3^k e^{\lambda_3(\xi)x}
\end{equation}
and we have also
\begin{equation}\label{eq:etilde}
\hat {\tilde e}^{1,k}(x,\xi)=\tilde\alpha_1^k e^{\lambda_1(\xi)x}\text{ and }\hat {\tilde e}^{2,k}(x,\xi)=\tilde\alpha_2^k e^{\lambda_2(\xi)x}+\tilde\alpha_3^k e^{\lambda_3(\xi)x}
\end{equation}
Using, $\hat{\cal G}\hat e^{1,1}=\hat{\cal G}\hat e^{2,1}$ we get:
\begin{equation}
\alpha_{\Gamma}:=\alpha_1^1(a(\xi)\bar c+\bar u R(\xi))=\alpha_2^1(a(\xi)\bar c-\bar uR(\xi))
\end{equation} 
Now we estimate the Fourier transform of the {\it correction} ${\tilde e}$. By using \eqref{eq:etilde} in the interface conditions of the {\it correction step} (\ref{step1}), we first get $\hat\gamma^1=i\xi\alpha_{\Gamma}R(\xi)$ and then:
\begin{equation}
\tilde\alpha_1^1=-\frac{\alpha_{\Gamma}}{a(\xi)+\bar u R(\xi)},\ \tilde\alpha_2^1=-\frac{\alpha_{\Gamma}}{a(\xi)\bar c-\bar uR(\xi)},\ \tilde\alpha_3^1=\frac{i\xi(\bar c^2-\bar u^2)}{a(\xi)}\cdot \frac{\alpha_{\Gamma}}{a(\xi)\bar c-\bar uR(\xi)}
\end{equation}
Now we estimate the Fourier transform of the {\it update} $e$. By using \eqref{eq:e} and \eqref{eq:etilde}   in the interface conditions of the {\it update step} (\ref{step2}), we first get that $\hat{\cal G}\hat e^{1,1}+\hat\delta^1=\hat{\cal G}\hat e^{2,1}+\hat\delta^1=0$ and $\hat {e}^{1,1}+\hat{\tilde e}^{1,1}=0$ on $\Gamma$. This means that the coefficients $(\alpha_i^2)_{1\le i\le 3}$ satisfy homogeneous systems of linearly independent equations. Therefore $\alpha_i^2=0$ for $i=1,2,3$ and then $\hat e^{1,2}=\hat e^{2,2}=0$. Therefore, the convergence is achieved in two steps.
\end{proof}
%
\section{A new algorithm applied to the Euler system}\label{schwarz-system}
After having found an optimal algorithm which converges in two steps for the third order model problem, we focus on the Euler system by translating this algorithm into an algorithm for the Euler system. It suffices to replace the operator ${\cal LG}$ by the Euler system and $Q$ by the last component  $F(W)_3$  of $F(W)$ in the boundary conditions. The algorithm reads:
\begin{algorithm}\label{algo1to2}
We choose the initial values $W^{1,0}$ and $W^{2,0}$ such that ${\cal G}F(W^{1,0})_3={\cal G}F(W^{2,0})_3$ and we compute $(W^{i,k+1})_{i=1,2}$ from $(W^{i,k})_{i=1,2}$ by the following iterative procedure:\\
{\bf Correction step} We compute the corrections $\tilde W^{1,k}$ and $\tilde W^{2,k}$ as solution of the homogeneous local problems:
\begin{equation}\label{step1to2}
\begin{array}{cc}
\left\{\begin{array}{l}
{\cal P}\tilde W^{1,k}=0\text{ in }\Omega_1,\\[2ex]
(A\nabla-\frac{1}{2}{\bf a}){\cal G}F(\tilde W^{1,k})_3\cdot{\bf n_1}=\gamma^k,\text{ on }\Gamma.
\end{array}\right.
&
\left\{\begin{array}{l}
{\cal P}\tilde W^{2,k}=0\text{ in }\Omega_2,\\[2ex]
(A\nabla-\frac{1}{2}{\bf a}){\cal G}F(\tilde W^{2,k})_3\cdot{\bf n_2}=\gamma^k,\text{ on }\Gamma,\\[2ex]
\tilde F(W^{2,k})_3=0,\text{ on }\Gamma.
\end{array}\right.
\end{array}
\end{equation} 
where $\gamma^k=-\frac{1}{2}\left[A\nabla{\cal G}F(W^{1,k})_3\cdot{\bf n_1}+A\nabla {\cal G}F(W^{2,k})_3\cdot{\bf n_2}\right]$.\\
{\bf Update step}.We update $W^{1,k+1}$ and $W^{2,k+1}$ by solving the local problems:
\begin{equation}\label{step2to2}
\begin{array}{cc}
\left\{\begin{array}{l}
{\cal P}W^{1,k+1}=f,\text{ in }\Omega_1,\\[2ex]
{\cal G}F(W^{1,k+1})_3={\cal G}F(W^{1,k})_3+\delta^k,\text{ on }\Gamma.
\end{array}\right.
&
\left\{\begin{array}{l}
{\cal P}\tilde W^{2,k+1}=f,\text{ in }\Omega_2,\\[2ex]
{\cal G}F(W^{2,k+1})_3={\cal G}F(W^{2,k})_3+\delta^k,\text{ on }\Gamma,\\[2ex]
F(W^{2,k+1})_3=F(W^{1,k})_3+F(\tilde W^{1,k})_3,\text{ on }\Gamma.
\end{array}\right.
\end{array}
\end{equation} 
where $\delta^k=\frac{1}{2}\left[{\cal G}F(\tilde W^{1,k})_3+{\cal G}F(\tilde W^{2,k})_3\right]$.
\end{algorithm}
This algorithm is quite complex since it involves second order derivatives of the unknowns in the boundary conditions on ${\cal G}F(W)_3$. It is possible to simplify it. By using the Euler equations in the subdomain, we have lowered the degree of the derivatives in the boundary conditions. After lengthy computations that we omit here, we find a simpler algorithm. We write it for a decomposition in two subdomains with an outflow velocity at the interface of domain $\Omega_1$ but with an interface not necessarily rectilinear. In this way, it is possible to figure out how to use for a general domain decomposition.

In the sequel, ${\bf n}=(n_x,n_y)$ denotes the outward normal to domain $\Omega_1$, $\partial_n=\nabla\cdot{\bf n}$ the normal derivative at the interface, $\partial_{\tau}=\nabla\cdot{\bf \tau}$ the tangential derivative, $U_n=U n_x+V n_y$ and $U_{\tau}=-U n_y+V n_x$ are respectively the normal and tangential velocity at the interface between the subdomains. Similarly, we denote $\bar u_n$ (resp. $\bar u_\tau$) the normal (resp. tangential) component of the velocity around which we have linearized the equations.
\begin{algorithm}\label{algo2}
We choose the initial values $W^{i,0}=(P^{i,0},U^{i,0},V^{i,0}),\,i=1,2$ such that $P^{1,0}= P^{2,0}$ and we compute $W^{i,k+1}$ from $W^{i,k}$ by the iterative procedure with two steps:\\
{\bf Correction step} We compute the corrections $\tilde W^{1,k}$ and $\tilde W^{2,k}$ as solution of the homogeneous local problems:
\begin{equation}\label{step21}
\begin{array}{cc}
\left\{\begin{array}{l}
{\cal P}\tilde W^{1,k}=0,\text{ in }\Omega_1,\\[2ex]
-(\beta +\bar u_\tau\partial_{\tau})\tilde U_{n}^{1,n}+\bar u_n\partial_{\tau}\tilde U_{\tau}^{1,k}=\gamma^k,\text{ on }\Gamma.
\end{array}\right.
&
\left\{\begin{array}{l}
{\cal P}\tilde W^{2,k}=0,\text{ in }\Omega_2,\\[2ex]
(\beta +\bar u_\tau\partial_{\tau})\tilde U_{n}^{2,k}-\bar u_n\partial_{\tau}\tilde U_{\tau}^{2,k}=\gamma^k,\text{ on }\Gamma\\[2ex]
\tilde P^{2,k}+\bar\rho\bar u_{n}\tilde U_{n}^{2,k}=0,\text{ on }\Gamma.
\end{array}\right.
\end{array}
\end{equation} 
where $\gamma^k=-\frac{1}{2}\left[(\beta +\bar u_\tau\partial_{\tau})(U_n^{2,k}-U_n^{1,k})+\bar u_n\partial_{\tau}(\tilde U_{\tau}^{1,k}-\tilde U_{\tau}^{2,k})\right]$.\\
{\bf Update step}.We compute the update of the solution $W^{1,k+1}$ and $W^{2,k+1}$ as solution of the local problems:
\begin{equation}\label{step22}
\begin{array}{cc}
\left\{\begin{array}{l}
{\cal P}W^{1,k+1}=f_1,\text{ in }\Omega_1,\\[2ex]
P^{1,k+1}=P^{1,k}+\delta^k,\text{ on }\Gamma.
\end{array}\right.
&
\left\{\begin{array}{l}
{\cal P}W^{2,k+1}=f_2,\text{ in }\Omega_2,\\[2ex]
P^{2,k+1}=P^{2,k}+\delta^k,\text{ on }\Gamma,\\[2ex]
(P+\bar\rho\bar u_{n}U_{n})^{2,k+1}=(P+\bar\rho\bar u_{n}U_{n})^{1,k}+(\tilde P +\bar\rho\bar u_{n}\tilde U_{n})^{1,k},\text{ on }\Gamma.
\end{array}\right.
\end{array}
\end{equation} 
where $\delta^k=\frac{1}{2}\left[\tilde P^{1,k}+\tilde P^{2,k}\right]$.
\end{algorithm}
\begin{proposition}\label{PropEquiv} 
For a domain $\Omega=\R^2$ divided into two non overlapping half planes, algorithms~\ref{algo1to2} and \ref{algo2} are equivalent and both converge in two iterations. 
\end{proposition}
\begin{proof} 
The Euler equations are invariant with respect to rotations so we can write at the interface the equations (\ref{euler3}) in the referential given by $({\bf n},{\bf\tau})$. We simply have to replace $\partial_x$(resp. $\partial_y$) by $\partial_n$(resp. $\partial_{\tau}$), $U$ (resp. $V$) by $U_n$ (resp. $U_\tau$) and $\bar u$ (resp. $\bar v$) by $\bar u_n$ (resp. $\bar u_\tau$). We denote by $\hat Q$ the Fourier transform along the interface of a function $Q$ and $\xi$ is the dual variable.\\
We first consider the boundary conditions \eqref{step2to2} of
Algorithm~\ref{algo1to2} and prove that they correspond to the ones
\eqref{step22}  of Algorithm~\ref{algo2}. First of all, from
\eqref{sec:F}, we have:
\begin{equation}
  \label{eq:FW3}
 \widehat{ F(W)}_3=-\frac{1}{\beta+i\xi \bar u_\tau} (\bar u_n \hat P + \bar \rho \bar u_n^2 \hat U_n)
\end{equation}
We apply now the operator ${\hat {\cal G}}$ and get:
\begin{equation}
  \label{eq:GFW3}
{\hat {\cal G}} (\widehat{ F(W)}_3)=-\frac{\bar u_n}{\beta+i\xi \bar u_\tau} ({\hat {\cal G}} (\hat P) + \bar \rho \bar u_n {\hat {\cal G}}(\hat U_n))
\end{equation}
By using the Euler equations satisfied by $W$ we can substitute ${\hat {\cal G}}(\hat U_n)$  with $-1/\bar \rho \partial_n \hat P$ (we can omit the right handside since it will appear on both sides of the boundary condition) and obtain:
\begin{equation}
  \label{eq:GFW3suite}
{\hat {\cal G}} (\widehat{ F(W)}_3)=-\frac{\bar u_n}{\beta+i\xi \bar u_\tau} ({\hat {\cal G}} (\hat P) -\bar u_n\partial_n \hat P)
= -\frac{\bar u_n}{\beta+i\xi \bar u_\tau} (i\bar u_\tau \xi + \beta )(\hat P)=-\bar u_n\hat P.
\end{equation}
Therefore, boundary condition of \eqref{step2to2} on the boundary of domain $\Omega_1$ reads:
\begin{equation}
  \label{eq:Prevealed}
  \bar u_n \hat P^{1,k+1}=\bar u_n\left(\hat P^{1,k}+\frac{1}{2}(\hat{\tilde P}^{1,k}+\hat{\tilde P}^{2,k})\right)
\end{equation}
By simplifying with $\bar u_n$, we get the boundary condition of \eqref{step22} on the pressure on the boundary of domain $\Omega_1$.\\ 
We now show how to obtain the second boundary condition in domain $\Omega_2$ both in the update step \eqref{step22} and in the correction step \eqref{step21}. From \eqref{eq:FW3}, we infer:
\begin{equation}
  \label{eq:Utaurevealed}
  \frac{1}{\beta+i\xi \bar u_\tau} (\bar u_n \hat{P}^{2,k+1} + \bar \rho \bar u_n^2 \hat{U}_n^{2,k+1})=\frac{1}{\beta+i\xi \bar u_\tau} (\bar u_n \hat{P}^{1,k} + \bar \rho \bar u_n^2 \hat{U}_n^{1,k}+ \bar u_n \hat{\tilde P}^{1,k} + \bar \rho \bar u_n^2 \hat{\tilde U}_n^{1,k})
\end{equation}
Multiplying by $(\beta+i\xi \bar u_\tau)/\bar u_n$, we  obtain the second boundary condition of \eqref{step21} and of  \eqref{step22}.\\
We now derive the first boundary condition of the correction step  \eqref{step21}  from the corresponding boundary condition of \eqref{step1to2}. We have to consider  
\begin{equation}
  \label{eq:robinFW3}
{\cal B}(W):=(A\nabla-\frac{1}{2}{\bf a}){\cal G}F(\tilde W)_3
\end{equation}
From \eqref{eq:Arond} and \eqref{eq:GFW3suite}, we
have:
\begin{equation}
  \label{eq:arondFW3}
  \widetilde{{\cal B}(W)}=-\bar u_n ((\bar c^2 -\bar u_n^2)\partial_n - \bar u_n (\beta+i \bar u_\tau \xi))(\hat P)
\end{equation}
In order to replace the normal derivative on $P$, we write the Euler system in the form:
\begin{equation}\label{eq:eulerx}
\partial_xW=-A^{-1}(\beta W+B\partial_y W-f)
\end{equation}
We get (once again omitting the right hand side $f$ that will appear
on both sides of the boundary conditions):
\begin{equation}
  \label{eq:dPdn}
  \partial_n P = \frac{\bar u_n}{\bar u_n^2-\bar c^2}[ -\bar u_n (\beta + \bar u_\tau \partial_\tau )(P) + \bar\rho\bar c^2(\beta+\bar u_\tau \partial_\tau)(U_n) -\bar u_n \bar\rho\bar c^2 \partial_\tau (U_\tau)]
\end{equation}
Using this equation in \eqref{eq:robinFW3}, 
\begin{equation}
  \label{eq:robinFW3fin}
  {\cal B}(W)=-\bar u_n \bar\rho\bar c^2 [(\beta+\bar u_\tau \partial_\tau)(U_n)-\bar u_n \partial_\tau (U_\tau)]
\end{equation}
To obtain  the first boundary condition of the correction step  \eqref{step21}, it suffices to multiply \eqref{step1to2} by $-\bar u_n \bar\rho\bar c^2$. This shows the equivalence of both algorithms. The convergence in two steps comes from the fact that Algorithm~\ref{algo1to2} was derived directly from Algorithm~\ref{algo1} which converges in two steps.
\end{proof}
\section{Discretization}\label{discrete}
In this section we will first present the discretization method used, a finite volume method on a uniform grid. Then we propose a strategy of discretization of the boundary conditions of the algorithm \ref{algo2} applied to the Euler system and we present some theoretical discrete estimates of the convergence of the method.
\subsection{A finite volume discretization}
We consider a domain $\Omega$ and the boundary value problem associated to (\ref{euler3}) with classical(natural) boundary conditions (see \cite{Dolean-etal:04:SYS}) on $\partial\Omega$. This  BVP  is discretized  using a   finite  volume scheme where   the flux  at  the interface  of the finite  volume cells  is  computed using a Roe type solver. We recall the method already described in \cite{Dolean:2004:OSM}. In order to  discretize  the BVP we consider  a regular quadrilateral grid where a  vertex $v_{ij}$  is characterized by
$$ 
v_{ij} = \left (\left(i - \frac{1}{2}\right)\Delta x,\left(j - \frac{1}{2}\right)\Delta y \right ),\,i,j\in\mathbb{Z}
$$
We associate to each vertex a finite  volume cell, $C_{ij} = \left [(i-1)\Delta x ~,~ i\Delta x\right ] \times \left [(j-1)\Delta x ~,~ j\Delta x\right ]$ which is a rectangle having as a center the vertex~ $v_{ij}$. A first order vertex centered finite volume formulation simply writes: 
\begin{equation}
\label{bilan}
\displaystyle\frac{W_{i,j}}{c\Delta t} + \frac{1}{|C_{ij}|} 
          \sum_{e\in\partial C_{ij}}|e|\Phi^e = f, 
\end{equation}
\noindent where $|C_{ij}|$ denotes the area of the cell $C_{ij}$, 
$|e|$ the length of the edge $e$ and $W_{i,j}$ the average value of the
unknown on the cell $C_{ij}$
$$ 
W_{i,j} = \displaystyle\frac{1}{|C_{ij}|}\int\limits_{C_{ij}} W(x,y)dxdy.
$$
Here, the elementary flux $\Phi_{ij}^e$ across edge $|e|$ is computed
by a Roe type scheme $\Phi^e = A_{{\bf n}}^+W_{i,j} + A_{{\bf n}}^-W_{k,l}$, where ${\bf n} = (n_x ~,~ n_y)$ is the outward normal to the the   edge $e$, $A_{{\bf n}}   =  n_xA_1 +  n_yA_2$   and $C_{kl}$  is the neighboring cell of $C_{ij}$  sharing the edge  $e$ with it. We can rewrite (\ref{bilan}) as
\begin{equation}
\label{bilan2}
\displaystyle\frac{W_{ij}}{c\Delta t} + 
\displaystyle\frac{|A_1|W_{i,j} + A_1^-W_{i+1,j} - A_1^+W_{i-1,j}}{\Delta x} + 
\displaystyle\frac{|A_2|W_{i,j} + A_2^-W_{i,j+1} - A_2^+W_{i,j-1}}{\Delta y} = 
f.
\end{equation}
 We will further denote $\bar \Delta x = \displaystyle\frac{\Delta x}{c\Delta t}$ and $\bar \Delta y = \displaystyle\frac{\Delta y}{c\Delta t}$, the non dimensioned counterpart of the mesh size in $x$ and $y$ directions.
\subsection{Discretizations of  the interface conditions}
In the following we consider the discretization of the interface conditions. In order to do that we will first write the semi-discrete system (only in the $y$ direction):
\begin{equation}\label{semidiscrete}
\left\{\begin{array}{l}
\beta \displaystyle P_{j}+\bar u\partial_x P_{j}+{\cal D}_{my}P_{j}+\bar\rho \bar c^2 \partial_x U_{j} + \bar\rho \bar c^2 {\cal D}_{py} V_{j}=f_1\\
\beta\displaystyle U_{j}+\bar u\partial_x U_{j}+ \bar v D^-_y U_{j}+\frac{1}{\bar\rho} \partial_x P_{j} =f_2\\
\beta\displaystyle V_{j}+\bar u \partial_x V_{j}+{\cal D}_{my}V_{j}+ \frac{1}{\bar\rho} {\cal D}_{py} P_{j} =f_3
\end{array}\right.
\end{equation}
where 
$$
\begin{array}{l}
{\cal D}_{my} = \displaystyle\left[\frac{\bar c+\bar
    v}{2}D^-_y-\frac{\bar c-\bar v}{2}D^+_y \right],\,{\cal D}_{py} =
\left[\frac{\bar c+\bar v}{2}D^-_y+\frac{\bar c-\bar v}{2}D^+_y \right]
\end{array}
$$
where we $M_n=\frac{\bar u}{\bar c}$ and $M_t=\frac{\bar v}{\bar c}$ are the normal and the tangential Mach number and $D^{\pm}_{y}$ are the usual finite difference operators. When writing the discrete counterpart of the algorithm \ref{algo2}, there is no need for a special treatment in the {\it update step} since the boundary conditions do not involve any derivatives. But one has to take care of the discretization of the interface conditions in the {\it correction step}. For this, we proceed as in the proof of proposition \ref{PropEquiv} but at the semi-discrete level (\ref{semidiscrete}). First of all, from (\ref{semidiscrete}) we obtain the discrete counterpart of (\ref{eq:dPdn}) 
\begin{equation}
\label{eq:dPdndis}
  \partial_n P = \frac{\bar u_n}{\bar u_n^2-\bar c^2}[ -\bar u_n
  (\beta + \bar u_\tau {\cal D}_{my} )(P) + \bar\rho\bar
  c^2(\beta+\bar u_\tau D_y^-)(U_n) -\bar u_n \bar\rho\bar c^2 {\cal D}_{py} (U_\tau)]
\end{equation}
and then by replacing it into (\ref{eq:arondFW3}) we obtain
\begin{equation}
\label{eq:robinFW3findis}
{\cal B}_{dis}(W)=\bar u_n \bar\rho\bar c^2 [-(\beta+\bar u_\tau
D_y^-)(U_n)+\bar u_n {\cal D}_{py}(U_\tau)+\bar u_{n}(\bar u_\tau\partial_{y}-{\cal D}_{my})\bar P].
\end{equation}
which suggests the use of the classical finite difference operator $D_y^-$ for the normal components of the velocity and ${\cal D}_{py}$ for the tangential one. Nevertheless, the question that remains open is which kind of approximation we should use for the pressure term. The most natural one is to cancel the pressure term. We will see in the discrete analysis that follows that it may lead to a diverging algorithm. Therefore, we will not completely drop this part but consider also the discretization of $(\bar u_\tau\partial_{y}-{\cal D}_{my})\bar P$ by $\Delta y D^+_yD^-_y\bar P$. As we shall see below, it leads to a converging algorithm. We say then that we have a stabilized algorithm. 

\paragraph{Theoretical discrete convergence results}

We will proceed to a discrete convergence analysis as in \cite{Dolean:2004:OSM} in order to decide which discretization of the boundary conditions is better. We will recall briefly the key ingredients of this analysis. We perform a discrete Fourier transform, by looking for the solution under the form:
\begin{equation}
\label{solu}
W_{i,j} = \sum_{\xi}\sum_{l=1}^3\alpha_{\xi l} 
          e^{(i-\frac{1}{2})\lambda_l(\xi)\Delta x}e^{Ij\xi\Delta y}V_l(\xi)
\end{equation} 
\noindent where $I^2=-1$. By   introducing     this  expression   into the   discrete   equation we get that for each  $\xi$, $\lambda_l(\xi)$ and $V_l(\xi)$
have to be the solution of
$$
\left (Id + \frac{|A_1| + A_1^{-}e^{\lambda_l(\xi)\Delta x} - 
                         A_1^+e^{-\lambda_l(\xi)\Delta x}}{\Delta x} + 
           \frac{|A_2| + A_2^{-}e^{I\xi\Delta y} - 
                         A_2^+e^{-I\xi\Delta y}}{\Delta y}\right ) 
V_l(\xi) = 0.
$$
If we denote by $L_l(\xi)=\frac{e^{-\lambda_l(\xi)}-1}{\Delta x}$ and by $e_y(\xi)= \frac{e^{I\xi\Delta y}-1}{\Delta x}$, we can further solve this systems numerically for each wave-number $\xi$. By  introducing these   expressions   in the interface conditions,  we get the discrete convergence rate.

\begin{figure}[htbp]
\begin{center}
\includegraphics[scale=0.25,angle=270]{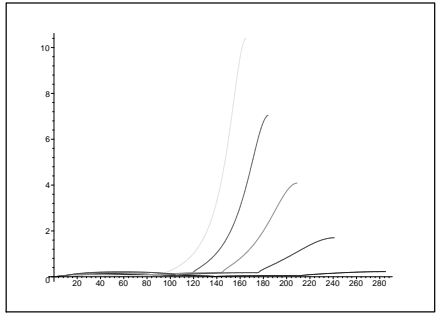}
\end{center}
\vspace{-2cm}
\caption{Convergence rate vs. Fourier number $\xi$ for different values of the normal Mach number, no stabilization used}
\label{FigStab1}
\end{figure}

\begin{figure}[htbp]
\begin{center}
\includegraphics[scale=0.25,angle=270]{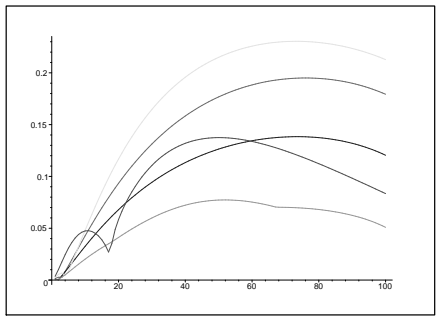}
\end{center}
\vspace{-2cm}
\caption{Convergence rate vs. Fourier number $\xi$ for different values of the normal Mach number, with stabilization term}
\label{FigStab2}
\end{figure}
In figures \ref{FigStab1} and \ref{FigStab2}, we show for a flow normal to the interface the convergence rate as a function of the wavenumber $\xi$ along the interface.  Two possibilities are considered: no stabilization and the stabilization mentioned above. These figures show the we need of the stabilization term for the high wavenumbers. 

\section{Numerical results}\label{results}
\begin{figure}[htbp]
\begin{center}
\includegraphics[scale=0.5]{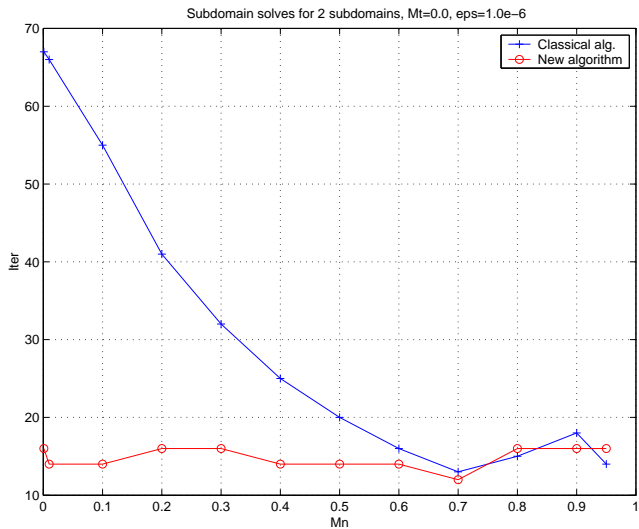}
\end{center}
\vspace{-0.5cm}
\caption{Subdomain solves counts for different values of $M_n$}
\label{comp1}
\end{figure}

\begin{figure}[htbp]
\begin{center}
\includegraphics[scale=0.5]{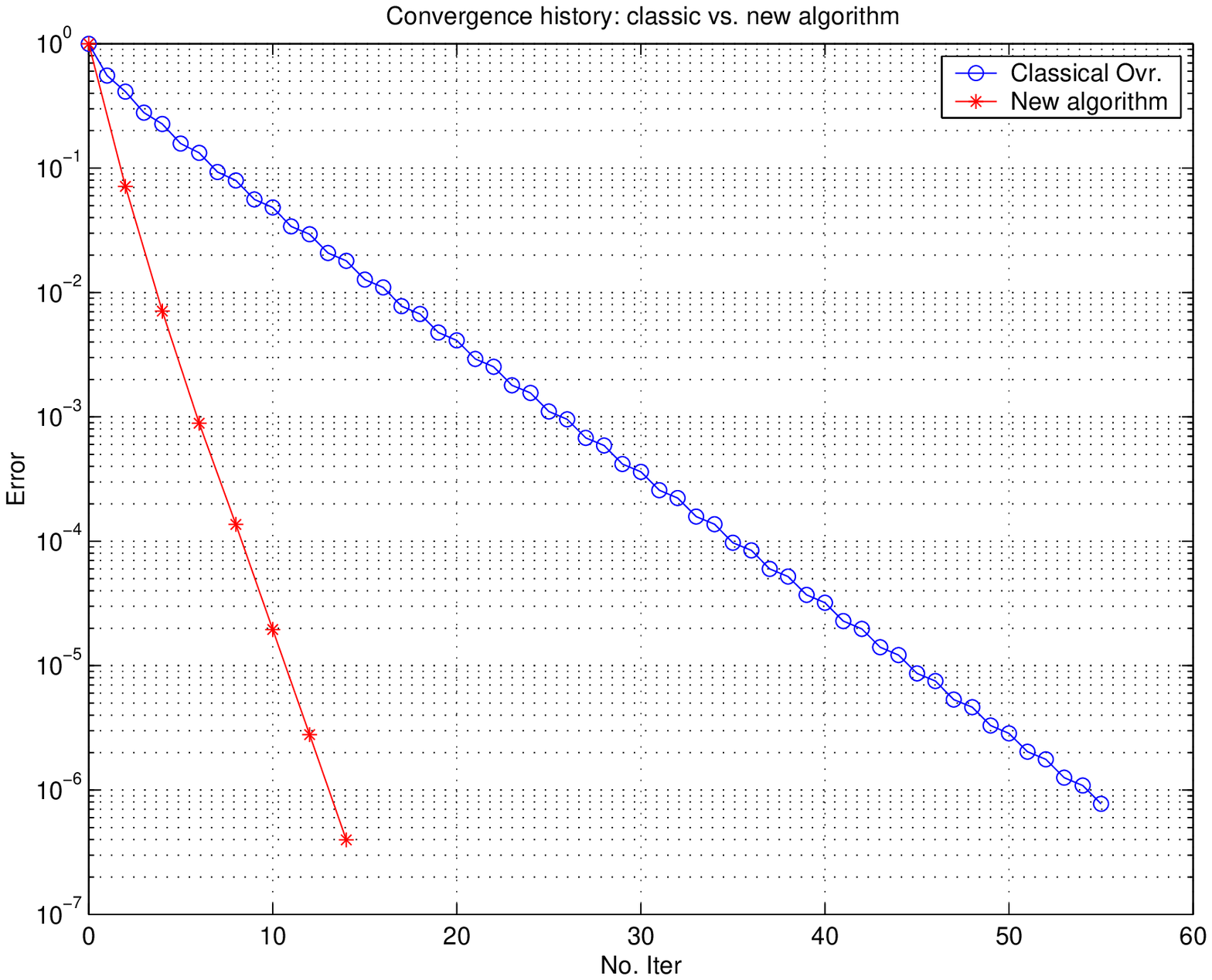}
\end{center}
\vspace{-0.5cm}
\caption{Convergence curves for the classical and the new algorithms}
\label{Hist1}
\end{figure}

\begin{figure}[htbp]
\begin{center}
\includegraphics[scale=0.5]{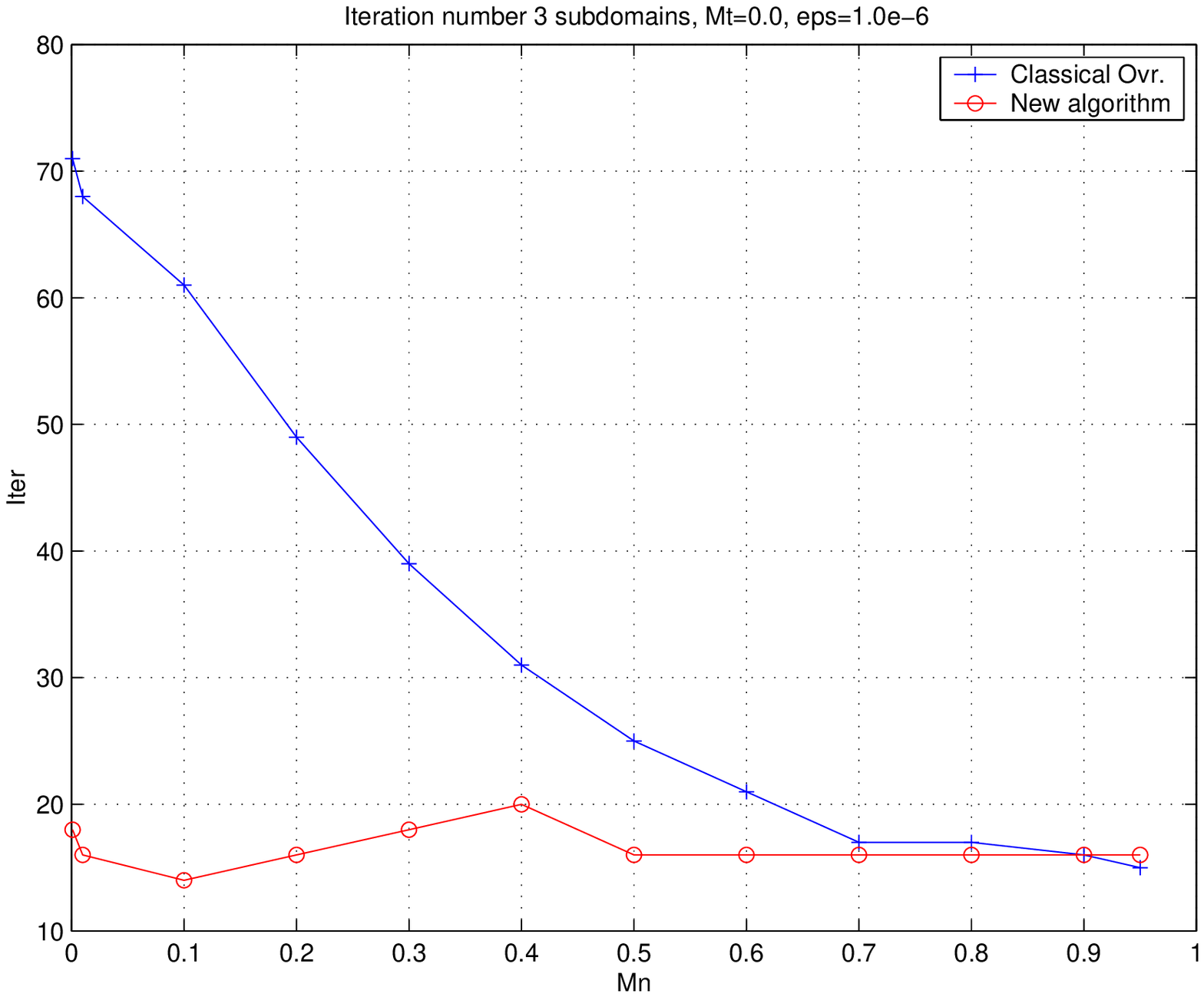}
\end{center}
\vspace{-0.5cm}
\caption{Convergence curves for the classical and the new algorithm}
\label{comp2}
\end{figure}

\begin{figure}[htbp]
\begin{center}
\includegraphics[scale=0.5]{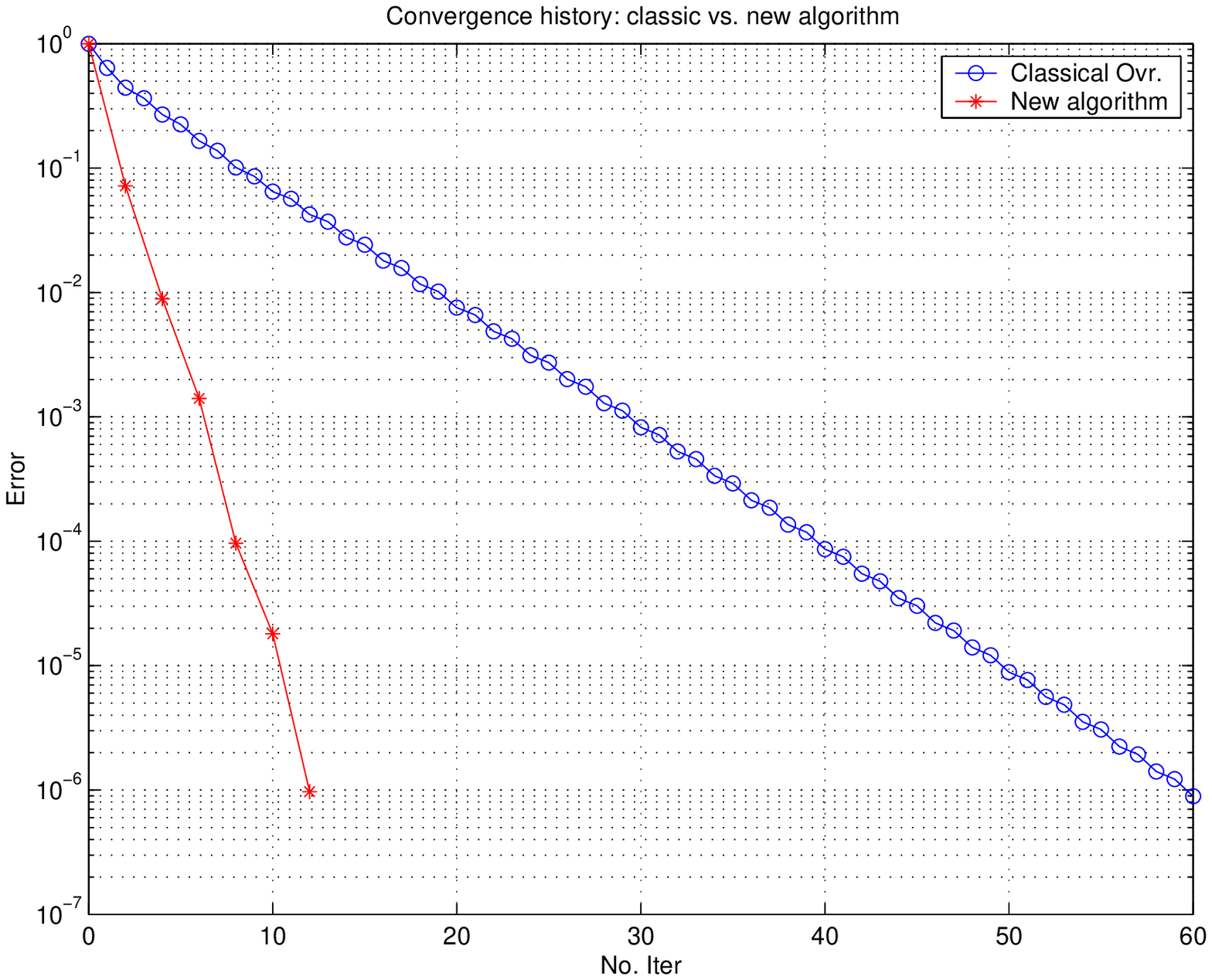}
\end{center}
\vspace{-0.5cm}
\caption{Convergence curves for the classical and the new algorithm}
\label{Hist2}
\end{figure}

We compare  the proposed method with the stabilization term and the classical method defined in \cite{Cai:1998:MORAS} or \cite{Dolean-etal:04:SYS} involving interface conditions   that  are  derived   naturally  from a  weak formulation of the  underlying boundary value  problem. We present here a set  of results of numerical experiments on a model problem. We consider a decomposition into different number of subdomains and for a linearization around a constant or non-constant flow. The computational  domain  is    given by   the   rectangle $[0~,~4]\times [0~,~1]$ with a uniform discretization using $80\times 20$ points. The numerical investigation  is limited to the solving of the linear system resulting from the first  implicit time step using a Courant number CFL=100. For all tests, the stopping criterion was a reduction of the maximum norm of the error by a factor $10^{-6}$. We consider first a decomposition into $2$ subdomains and a linearization around a constant state and we are solving the homogeneous equations satisfied by the error vector at the first time step. In the following, for the new algorithm, each iteration counts for $2$ as we need to solve twice as much local problems than with the classical algorithm. For an easier comparison of the algorithms, the figures shown in the tables are the number of subdomains solves. Table \ref{Tab1} summarizes the number of subdomains solves for different values of the  reference Mach number for the new and the classical algorithm for a normal flow to the interface $M_t=0.0$. The same results are presented in figure \ref{comp1}. 
\begin{table}
\begin{center}
\begin{tabular}{|c|c|c|c|c|c|}
\hline 
$M_n $ & Classical & New DDM & $M_n $ & Classical & New DDM  \\
\hline 
 0.001 & 67 & 18 &  0.4  & 25 & 14\\
 0.01  & 66 & 16 &  0.5  & 20 & 14\\ 
 0.1   & 55 & 14 &  0.6  & 16 & 14\\
 0.2   & 41 & 16 &  0.7  & 13 & 12\\
 0.3   & 32 & 16 &  0.8  & 15 & 14\\
\hline 
\end{tabular} 
\caption{Subdomain solves counts for different values of $M_n$} \label{Tab1}
\end{center} 
\end{table} 
In Table \ref{Tab2}, we consider a linearization around a variable state where the tangential velocity is given by $M_t(y) = 0.1(1+\cos(\pi y))$ and we vary the normal Mach number.
\begin{table}
\begin{center}
\begin{tabular}{|c|c|c|c|c|c|}
\hline 
$M_n $ & Classical & New DDM & $M_n $ & Classical & New DDM  \\
\hline 
 0.001 & 32 & 16 &  0.4  & 18 & 14\\
 0.01  & 30 & 16 &  0.5  & 16 & 12\\ 
 0.1   & 28 & 14 &  0.6  & 15 & 12\\
 0.2   & 24 & 14 &  0.7  & 14 & 12\\
 0.3   & 20 & 14 &  0.8  & 14 & 14\\
\hline 
\end{tabular}
\caption{Subdomain solves counts for different values of $M_n$, $M_t(y)$} \label{Tab2} 
\end{center} 
\end{table} 
In figure \ref{Hist1}, we linearize the equations around a variable state for a normal flow to the interface ($M_t=0.0$), where the initial normal velocity is gives by $M_n(y) = 0.5(0.2+0.04\tanh(y/0.2)))$. 
The sensitivity to the mesh size is shown in the Table \ref{Tab21}. We can see that for the new algorithm the growth in the number of iterations is very weak as the mesh is refined, the same property being already known for the classical one.
\begin{table}
\begin{center}
\begin{tabular}{|c|c|c|c|c|c|}
\hline 
$h$ ($M_n=0.001$) & Classical & New DDM & $h$ ($M_n=0.1$) & Classical & New DDM  \\
\hline 
 1/10 & 65 & 18 &  1/10  & 56 & 12\\
 1/20 & 67 & 18 &  1/20  & 57 & 14\\ 
 1/40 & 70 & 18 &  1/40  & 59 & 16\\
\hline 
\end{tabular}
\caption{Subdomain solves counts for different mesh size} \label{Tab21} 
\end{center} 
\end{table} 
The next tests concern a decomposition into $3$ subdomains. We first consider a linearization around a constant state. Table \ref{Tab3} gives the results  for different values of the  reference Mach number for the new and the classical algorithm. The same results are presented in figure \ref{comp2}.
\begin{table}
\caption{Classical vs. the new algorithm}  
\begin{center}
\begin{tabular}{|c|c|c|c|c|c|}
\hline 
$M_n $ & Classical & New DDM & $M_n $ & Classical & New DDM  \\
\hline 
 0.001 & 71 & 18 &  0.4  & 31 & 20\\
 0.01  & 68 & 16 &  0.5  & 25 & 16\\ 
 0.1   & 61 & 14 &  0.6  & 21 & 16\\
 0.2   & 49 & 16 &  0.7  & 17 & 16\\
 0.3   & 39 & 18 &  0.8  & 16 & 16\\

\hline 
\end{tabular} 
\caption{Subdomain solves counts for different values of $M_n$} \label{Tab3}
\end{center} 
\end{table}
In figure \ref{Hist2}, we consider for a three subdomains decomposition a linearization around a variable state for a normal flow to the interface. The normal velocity is given by $M_n(y) = 0.5(0.2+0.04\tanh(y/0.2)))$ (the same as for the $2$ subdomain case).
These tests show that the new algorithm is very stable with respect to various parameters such as the mesh size and the Mach number. We see that the convergence in two iterations of the continuous algorithm is lost at the discrete level although the subdomain solves are very reasonable. Moreover, a stabilization was necessary for the discretization of the interface condition \eqref{step21} in order to keep the algorithm converging. The optimal discretization of this interface condition is not yet quite well known. The comparison with the classical algorithm is favorable for  Mach numbers smaller than $0.5$ and especially very low Mach numbers by a factor of almost $4$.

\section{Conclusion}
In this paper we designed a new domain decomposition for the Euler equations inspired by the idea of the Robin-Robin preconditioner \cite{Achdou:1997:RPA} applied to the advection-diffusion equation. We used the same principle after reducing the system to scalar equations via a Smith factorization. The resulting algorithm behaves very well for low Mach numbers, where usually the classical algorithm doesn't give very good results. We reduce the number of subdomain solves by almost a factor $4$ for linearization around a constant and variable state as well. A general theoretical study and more comprehensive numerical tests have to be done in order to firmly assess the applicability of the proposed algorithm to large scale computations.\\ This work can also be seen as a first step for deriving new domain decomposition methods for the 2D or 3D compressible Navier-Stokes equations. Indeed, the derivation of the algorithm which is based on the Smith factorization is in fact general and can be applied to arbitrary systems of partial differential equations, see \cite{Nataf:2005:NDS}  for the Stokes and Oseen equations and \cite{Dolean:2005:NCD}  for a general presentation. 

\addcontentsline{toc}{section}{References}
\bibliography{newddm,ddm,frednewddm}
\end{document}